\documentclass[a4paper,12pt]{article}
\usepackage{tabularx} % extra features for tabular environment
\usepackage{amsmath,amsthm,verbatim,amssymb,amsfonts,amsopn,amscd}
\usepackage{graphicx} % takes care of graphic including machinery
\usepackage[margin=1in,a4paper]{geometry} % decreases margins

\usepackage[square,numbers]{natbib}
\usepackage[utf8]{inputenc}
\usepackage[T1]{fontenc}
\usepackage[UKenglish]{babel}
\usepackage{enumerate}
\usepackage{paralist} % for asparaenum - numbering as paragraphs in proofs
\usepackage[final]{hyperref} % adds hyper links inside the generated pdf file
\hypersetup{
	colorlinks=true,       % false: boxed links; true: colored links
	linkcolor=blue,        % color of internal links
	citecolor=blue,        % color of links to bibliography
	filecolor=magenta,     % color of file links
	urlcolor=blue         
}

\usepackage{tikz-cd}
\tikzcdset{row sep/normal=2.2em}
\tikzcdset{arrow style=math font}

\theoremstyle{plain}
\newtheorem{theorem}{Theorem}[section]
\newtheorem{corollary}[theorem]{Corollary}

\newtheorem{proposition}[theorem]{Proposition}

\theoremstyle{definition}
\newtheorem{definition}[theorem]{Definition}

\DeclareMathOperator{\dom}{dom}

\DeclareMathOperator{\Sub}{Sub}

\DeclareMathOperator{\ev}{ev} %Evaluation
\DeclareMathOperator{\pr}{pr}

\newcommand{\set}[2]{\{{#1}\,\mid\,{#2}\}}

\DeclareMathOperator{\calA}{\mathcal{A}}

\DeclareMathOperator{\calS}{\mathcal{S}}

\newcommand{\IN}{\mathbb{N}}
\newcommand{\IZ}{\mathbb{Z}}

\newcommand{\IR}{\mathbb{R}}

\newcommand{\IT}{\mathbb{T}}

\newcommand{\into}{\hookrightarrow}
\newcommand{\monto}{\rightarrowtail}

\DeclareMathOperator{\iAlg}{\mathbf{-iAlg}}

\DeclareMathOperator{\Alg}{\mathbf{-Alg}}
\DeclareMathOperator{\Mod}{\mathbf{-Mod}}

\DeclareMathOperator{\Set}{\mathbf{Set}}

\DeclareMathOperator{\IIT}{I[\mathbb{T}]}

\DeclareMathOperator{\End}{End}

\renewcommand{\vec}[1]{\overrightarrow{#1}}
\newcommand{\la}{\langle}
\newcommand{\ra}{\rangle}
\title{Gluing of infinitesimal models of algebraic theories}
\author{\href{mailto:filip.bar.research@gmail.com}{Filip B\'{a}r}}
\date{\itshape In memoriam of my mother (22.04.1955 - 06.12.2022)}
\begin{document}
  \maketitle

  \begin{abstract}
    \noindent 
    Categories of models of algebraic theories have good categorical properties except for gluing. Building upon insights and examples from Synthetic Differential Geometry, we introduce a generalisation of models of algebraic theories to infinitesimal models. We demonstrate that the category of infinitesimal models retains most of the good categorical properties, but with a stark improvement in the behaviour of gluing. This makes infinitesimal models an interesting natural construction with the ability to interpolate between algebra and geometry.
  \end{abstract}
  
  \section{Introduction}
  In \cite{Kock:Levi-Civita} Kock has shown that a (formal) manifold in Synthetic Differential Geometry admits affine combinations of points that are pairwise mutual infinitesimal neighbours. In this and subsequent work \cite{Kock:Synthetic_Geometry_Manifolds, Kock:flatness} he has made an extensive use of this geometric algebra of infinitesimally affine combinations linking it with well-known concepts and constructions from Differential Geometry. Building on Kock's work the author has been trying to understand in which sense formal manifolds are models of affine spaces and whether this can be extended to other algebraic theories like groups and vector spaces. This has led him to formulate the notion of an infinitesimal model of an algebraic theory as a space equipped with an infinitesimal structure that serves as the domain of the operations of the theory in \cite{Bar:thesis}. 

  Similar to the structure of a topology an infinitesimal structure defines the collections of infinitesimally neighbouring points. However, in contrast to open sets the collections are finite tuples of points. This allows to use an infinitesimal structure as a domain for operations of an algebraic theory\footnote{In this paper we shall assume that all algebraic theories are finitary, i.e. all operations have a finite arity. Moreover, we shall take the "working" definition of an algebraic theory via its presentation, i.e. a signature with one sort, function symbols and a set of equations.}. In \cite{Bar:second_order_affine_structures} the author has given a first major application showing symmetric affine connections equivalent to second-order infinitesimally affine structures on a manifold, making full use of infinitesimally affine spaces as additional structures that go beyond the above mentioned property of formal manifolds admitting affine combinations of mutual (first-order) infinitesimal neighbours utilised by Kock. In \cite{Bar:2nd-order_i-groups} this equivalence could be extended to non-symmetric affine connections, linking the latter with infinitesimal models of groups at every point of the manifold. The many examples of infinitesimal models of affine spaces, vector spaces, groups and Lie algebras encountered in Synthetic Differential Geometry, as well as their applications, justify a general theory of infinitesimal algebra.
  
  The aim of this paper is to introduce the notion of an infinitesimal model of an algebraic theory $\IT$ as a natural construct that subsumes and extends the notion of a $\IT$-model in its own right and independent of geometry; the hope being that it finds applications in other fields besides Synthetic Differential Geometry. To serve this purpose we investigate the properties of categories of infinitesimal models. By relying on some strong general results from Categorical Logic we will be able to quickly establish the general categorical properties. The main focus will lie on the particular gluing properties of infinitesimal models, i.e. which colimits of infinitesimal models can be computed from the underlying sets. We shall see that infinitesimal models have some remarkable gluing properties: in stark contrast to $\IT$-models coproducts become essentially unions of the underlying sets, while coequalizers are, in general, constructed from quotients of congruences like in the category of $\IT$-models. Which other colimits are lifted depends on how `big' the infinitesimal structures of the models are in relation to each other. We shall provide sufficient conditions in our first and third gluing theorem and their corollaries. This shows how infinitesimal models interpolate between algebra and geometry.    
  
  In line with the aim of this paper we will use the category $\Set$ of sets as the base category. For the purpose of Synthetic Differential Geometry, however, we need to work over a well-adapted model \cite{Dubuc:Well-adapted_models} as a base, which is a Grothendieck topos. Indeed, all the results stated in this paper generalise to a Grothendieck topos. The proofs of the gluing theorems and their corollaries transfer by re-phrasing them within the internal language of the topos.       
 
%%%%%%%%%%%%%%%%%%%%%%%%%%%%%%%%%%%%%%%%%%%%%%%%%%%%%%%%%%%%%%%%%%%%%%%%%%%%%%%%%%%%%%%%%%%%%%%%%%%%%5

  \section{Infinitesimal structures}

  An \emph{infinitesimal structure} on a set $A$ (or `i-structure' for short) amounts to give an $n$-ary relation $A\la n\ra$ for each $n\in\IN$ that defines which $n$ points in $A$ are considered as being `infinitesimally close' to each other.
	\begin{definition}[i-structure]\label{def:i-structure}
		Let $A$ be a space. An \emph{\textbf{i-structure}} on $A$ is an $\IN$-indexed family $n\mapsto A\la n\ra\subseteq A^n$ such that
		\begin{enumerate}[(1)]
			\item $A\la 1\ra = A$, $A\la 0\ra=A^0=1$ (the `one point' space, or terminal object)
			\item For every map $h:m\to n$ of finite sets 
			and every $(P_1,\ldots,P_n)\in A\la n\ra$ we have $ (P_{h(1)},\ldots,P_{h(m)})\in A\la m\ra$ 
		\end{enumerate}
	\end{definition}
	The first condition is a normalisation condition. The second condition makes sure that the relations are compatible: if we have a family of points that are infinitesimally close to each other, then so is any subfamily of these points, or any family created from repetitions. In particular, we obtain that the $A\la n\ra$ are symmetric and reflexive relations. An $n$-tuple $(P_1,\ldots,P_n)\in A^n$ that lies in $A\la n\ra$ will be denoted by $\la P_1,\ldots,P_n\ra$. A map $f:A\to X$ that preserves infinitesimal structures, i.e. satisfies $f^n(A\la n\ra)\subseteq X\la n\ra$, is called an \emph{i-morphism}.
	
	Two trivial examples of i-structures on $A$ are the discrete and the indiscrete i-structure obtained by taking $A\la n\ra$ to be the diagonal $\Delta_n$, respectively the whole $A^n$. A non-trivial class of examples, and also the i-structures that are of main interest in Synthetic Differential Geometry are the i-structures generated by the first neighbourhood of the diagonal (as relations). For example, let $R$ be a commutative ring. Recall that
	$$
		D(n)=\set{(d_1,\ldots,d_n)\in R^n}{d_id_j=0,\ 1\leq i,j\leq n}
	$$
	On $R^n$ the first neighbourhood of the diagonal is given by 
	$$
	\set{(P_1,P_2)}{P_2-P_1\in D(n)}
	$$ 
	This is a symmetric and reflexive relation and we can construct an i-structure from it: take the first neighbourhood of the diagonal as $R^n\la 2\ra$ and define the \emph{nil-square i-structure} on $R^n$ by
	$$
		R^n\la m\ra=\set{(P_1,\ldots,P_m)}{(P_i,P_j)\in R^n\la 2\ra,\ 1\leq i,j \leq m}
	$$    
	This i-structure is thus \emph{generated} by $R^n\la 2\ra$. Not all i-structures $A\la-\ra$ of interest need to be generated by $A\la 2\ra$. The first- and second-order i-structures defined in \cite{Bar:2nd-order_i-groups, Bar:second_order_affine_structures} are not, for example. In \cite[prop.~17.4]{Kock:SDG} Kock has shown that every formal manifold carries a natural nil-square i-structure glued together from its local models in $R^n$, where $R$ satisfies the Kock-Lawvere axiom scheme.\footnote{This can also be deduced from our third gluing theorem~\ref{thm:gluing_local_models} when a well-adapted model is used as the base category instead of $\Set$.}
 
 %%%%%%%%%%%%%%%%%%%%%%%%%%%%%%%%%%%%%%%%%%%%%%%%%%%%%%%%%%%%%%%%%%%%%%%%%%%%%%%%%%%%%%%%%%%%%%%%%%%%%5
  
  \section{Clones and algebraic theories}

  The idea behind an infinitesimal model $A$ of an algebraic theory $\mathbb{T}$ is that, rather than on the products of the underlying space, any $n$-ary operation is only defined on $A\la n\ra$ for a given i-structure on $A$. To be able to define this formally we require a representation of an algebraic theory that considers all the operations of a theory $\mathbb{T}$ (including the derived ones) sorted by arity. We shall use the structure of an abstract \emph{clone} for this purpose. Our definition is based on \cite[def.~1.2.1]{Gould:PhD}:

  \begin{definition}[\emph{Clone}]\label{def:clone}
    The data of a \textbf{clone} $O$ consists of:
    \begin{itemize}
      \item For every $n\in \IN$ a set $O(n)$. 
      \item For every $(n,k)\in\IN^2$ a map $\ast_{nk}:O(n)\times O(k)^n \to O(k)$. 
      \item For every $n\geq 1$ and $1\leq j\leq n$ elements $\pi^n_j\in O(n)$. 
    \end{itemize}
   
    \begin{enumerate}[(1)]
      \item (\textbf{Associativity})\hspace{0.3em} For every $\sigma\in O(n)$, $t_1,\ldots,t_n\in O(m)$, $s_1,\ldots s_m\in O(k)$
      \begin{align*}
        \sigma\ast_{nk}(t_1\ast_{mk}(s_1,\ldots, s_m),\ldots,t_n\ast_{mk}(s_1,\ldots, s_m))
        =(\sigma\ast_{nm}(t_1,\ldots,t_n))\ast_{mk}(s_1,\ldots,s_m).
      \end{align*} 
      In particular, the naming and evaluation of constants are compatible as well as the naming of constant operations: 
      \begin{itemize}
        \item In the case $m=0$ the $t_i$ are constants and the associativity states
        $$\sigma\ast_{nk}(t_1,\ldots,t_n)
        =\ast_{0k}(\sigma\ast_{n0}(t_1,\ldots,t_n))$$
        \item In the case $n=0$ the operation $\sigma$ is a constant and the associativity states
        $$\ast_{0k}(\sigma) = (\ast_{0m}(\sigma))\ast_{mk}(s_1,\ldots,s_m)$$
      \end{itemize}
      \item (\textbf{Projection})\hspace{0.3em} For every $n\geq 1$, $1\leq j\leq n$, $t_1,\ldots,t_n\in O(m)$
      $$ \pi^n_j\ast_{nm}(t_1,\ldots,t_n)=t_j$$
      
      \item (\textbf{Unit})\hspace{0.3em} For each $\sigma\in O(n)$, $n\geq 1$
      $$\sigma\ast_{nn}(\pi^n_1,\ldots,\pi^n_n)=\sigma$$ 
      
      \item (\textbf{Normalisation}) \hspace{0.3em} $\ast_{00}(c)=c.$
    \end{enumerate}
  \end{definition}
  A \emph{clone homomorphism} $f:O\to O'$ is a family of functions $f_n:O(n)\to O'(n)$ for each $n\in\mathbb{N}$ commuting with the clone operations. A \emph{clone algebra} is an action of a clone $O$ on a set $A$: 

  \begin{definition}[\emph{Clone algebra}]\label{def:clone_algebra}
    Let $O$ be a clone. An \textbf{$O$-algebra} $A$ consists of:
    \begin{itemize}
      \item A set $A$.
      \item For every $n\in\mathbb{N}$ a map $\bullet_{n}:O(n)\times A^n \to A$.  
    \end{itemize}
    satisfying
    \begin{enumerate}[(1)]
      \item (\textbf{Associativity})\hspace{0.3em} For every $\sigma\in O(n)$, $t_1,\ldots,t_n\in O(m)$, $a_1,\ldots a_m\in A$
      \begin{align*}
        \sigma\bullet_{n}(t_1\bullet_{m}(a_1,\ldots, a_m),\ldots,t_n\bullet_{m}(a_1,\ldots, a_m))
        =(\sigma\ast_{nm}(t_1,\ldots,t_n))\bullet_{m}(a_1,\ldots,a_m).
      \end{align*} 
      In particular, the action preserves constants: 
        $$\bullet_{0}(c) = (\ast_{0m}(c))\bullet_{m}(a_1,\ldots,a_m)$$
      \item (\textbf{Projection})\hspace{0.3em} For every $n\geq 1$, $1\leq j\leq n$, $a_1,\ldots,a_n\in A$
      $$ \pi^n_j\bullet_{n}(a_1,\ldots,a_n)=a_j$$
    
    \end{enumerate}
  \end{definition}
  
  As it is familiar from representation theory, an $O$-algebra structure on $A$ is equivalent to a clone homomorphism $ O\to\End(A)$, where $\End(A)$ denotes the endomorphism clone defined by (multi-) composition of maps $A\times\ldots\times A\to A$. 

  \emph{$O$-algebra homorphisms} are the `equivariant' maps $f:A\to A'$, i.e. maps between the underlying sets that commute with the respective actions; we obtain a category $O\Alg$ of $O$-algebras and $O$-algebra homomorphisms.

  The (abstract) clone\footnote{One can also form the clone of operations of $\mathbb{T}$ for a $\mathbb{T}$-model $A$. This is how clones have been introduced in universal algebra, originally \cite[chap.~III.3]{Cohn:Univ_Alg}. However, we are interested in the clone encoding $\mathbb{T}$ rather than just one of its models.} $O_{\mathbb{T}}$ corresponding to an algebraic theory $\mathbb{T}$ can be obtained as follows: Given a presentation $(\Sigma,E)$ of $\mathbb{T}$ define $O_{\mathbb{T}}(n)$ as the finitely generated free $\mathbb{T}$-algebra $F_{\mathbb{T}}(n)=T_\Sigma(n)/E_n$, where the $T_\Sigma(n)$ is the term algebra of terms over the signature $\Sigma$ in $n$ variables. The operation of substitution induces maps
  $\ast_{nm}:T_\Sigma(n)\times T_\Sigma(m)^n\to T_\Sigma(m)$ compatible with the congruence relations and thus descends to a map
  $$\ast_{nm}: F_\mathbb{T}(n)\times F_{\mathbb{T}}(m)^n\to F_{\mathbb{T}}(m),\quad ([t],([s_1],\ldots,[s_n])) \mapsto [t[s_1/x_1,\ldots,s_n/x_n]]$$
  The $\pi^n_j$ name the (equivalence classes of) variables $[x_j]\in F_{\mathbb{T}}(n)$ for $n\geq 1$. The axioms of a clone follow from the corresponding properties of substitution of terms. Regarding normalisation in particular, since $F_{\mathbb{T}}(0)$ is the set of (equivalence classes of) terms with no free variables, substitution becomes the identity map.

  Conversely, given a clone $O$ we can can construct an algebraic theory $\mathbb{T}_O$ by defining an $n$-ary function symbol for every operation $\sigma\in O(n)$, $n\in\mathbb{N}$ and take as the set of equations all the defining equations of $O$ stated in definition~\ref{def:clone}. These constructions induce isomorphisms between the concrete categories of models of the theory and the algebras of the corresponding clone:

  \begin{proposition}[Clones and algebraic theories]\label{prop:compare_syntactic_clone}
    Clones and algebraic theories are equivalent in the following sense:
  \begin{enumerate}[(1)]
    \item For every algebraic theory $\mathbb{T}$ there is a clone $O_\mathbb{T}$ such that $\mathbb{T}\Mod$ and $O_\mathbb{T}\Alg$ are isomorphic categories over $\Set$.
    \item For every clone $O$ there is an algebraic theory $\mathbb{T}_O$ such that $\mathbb{T}_O\Mod$ and $O\Alg$ are isomorphic categories over $\Set$.
    \item The clones $O_{\mathbb{T}_O}$ and $O$ are isomorphic.
  \end{enumerate}
  \end{proposition}
  \begin{proof}
    See \cite[thm.~1.4.1]{Bar:thesis} for the remaining parts that need to be shown.
  \end{proof}

%%%%%%%%%%%%%%%%%%%%%%%%%%%%%%%%%%%%%%%%%%%%%%%%%%%%%%%%%%%%%%%%%%%%%%%%%%%%%%%%%%%%%%%%%%%%%%%%%%%%%5

  \section{Infinitesimal models of algebraic theories}

  We are now ready to define an infinitesimal model of an algebraic theory $\mathbb{T}$ as an infinitesimal algebra of the clone $O_\mathbb{T}$ that acts on an i-structure $A$. Apart from the domain of the operations the main difference to (total) $O$-algebras is the neighbourhood axiom, which guarantees that operations on infinitesimally neighbouring points result in infinitesimally neighbouring points again. This is necessary to be able to define associativity.

  \begin{definition}[\emph{i-algebra of a clone}]\label{def:i-alg}
    Let $O$ be a clone. An i-structure $A$ together with a family of maps
    $$\bullet_n: O(n)\times A\la n\ra \to A\la 1\ra,\qquad n\in\IN$$
    is an \emph{\textbf{infinitesimal $O$-algebra}} if it satisfies the following axioms:
    \begin{enumerate}[(1)]
	    \item (\textbf{Neighbourhood})\hspace{0.3em} For each pair $(n,m)\in\IN^2$, $n\geq 1$, $\sigma_1,\ldots, \sigma_n\in O(m)$, and $a\in A\la m\ra$ we have 
        $$
          \la\sigma_1\bullet_m a,\ldots,\sigma_n\bullet_m a\ra\in A\la n\ra
        $$
	      Note that in the case of constants ($m=0$) this becomes $\la\bullet_0 (\sigma_1),\ldots,\bullet_0(\sigma_n)\ra\in A\la n\ra$.
	
	    \item (\textbf{Associativity})\hspace{0.3em} For each pair $(n,m)\in\IN^2$, $\sigma\in O(n)$, $t_1,\ldots,t_n\in O(m)$, $\la a_1,\ldots,a_m\ra\in A\la m\ra$
      \begin{align*}
        \sigma\bullet_{n}(t_1\bullet_{m}(a_1,\ldots, a_m),\ldots,t_n\bullet_{m}(a_1,\ldots, a_m))
        =(\sigma\ast_{nm}(t_1,\ldots,t_n))\bullet_{m}(a_1,\ldots,a_m).
      \end{align*} 
    
      \item (\textbf{Projection})\hspace{0.3em} For every $n\geq 1$, $1\leq j\leq n$, $\la a_1,\ldots,a_n\ra\in A\la n\ra$
      $$ \pi^n_j\bullet_{n}(a_1,\ldots,a_n)=a_j$$ 
    \end{enumerate}
  \end{definition}

  An i-$O$-homomorphism $h:(A,\bullet)\to (A',\bullet')$ is an i-morphism $h:A \to A'$ that commutes with the operations, i.e. 
  $$ h(\sigma\bullet_n x)=\sigma\bullet'_n h^n(x),\qquad x\in A\la n\ra$$
  Infinitesimal $O$-algebras and infinitesimal $O$-homomorphisms form a category $O\iAlg$. Due to the indiscrete i-structure every $O$-algebra is also an infinitesimal $O$-algebra, which we shall refer to as a \emph{total (i-)$O$-algebra}; the category $O\Alg$ is thus a full subcategory of $O\iAlg$ and infinitesimal models extend the notion of $\mathbb{T}$-models. 

  Interesting examples of infinitesimal models of algebraic theories arise naturally in Synthetic Differential Geometry and in Algebraic Geometry over rings and $C^\infty$-rings\footnote{In each of these examples we need to work over base categories that are different from $\Set$ though.}. For example, any formal manifold $M$ \cite[chap.~17]{Kock:SDG} is an infinitesimal model of the theory of affine combinations over a commutative $\mathbb{R}$-algebra $R$ satisfying the Kock-Lawvere axiom scheme \cite[thm.~3.2.8]{Bar:thesis}\footnote{This result is foreshadowed by \cite[thm.~2.2]{Kock:Levi-Civita}. Note that the proof of \cite[thm.~3.2.8]{Bar:thesis} relies on \cite[thm.~2.6.19]{Bar:thesis} for which we provide a counterexample in this paper; the third gluing theorem~\ref{thm:gluing_local_models} should be used instead.}. For any point $P\in M$ the subspace of (first-order) infinitesimal neighbours of $P$ carries the structure of an infinitesimal $R$-module. If $M$ is also a group, then the subspace of infinitesimal neighbours of the neutral element $e$ forms an infinitesimal group (but not for the nil-square i-structure) \cite[thm.~3.4]{Bar:2nd-order_i-groups}. To make these examples more relatable note that any well-adapted model of Synthetic Differential Geometry is a Grothedieck topos equipped with a fully faithful embedding of the category of $C^\infty$-manifolds. Moreover, it maps $\mathbb{R}$ to $R$ and manifolds to formal manifolds \cite{Dubuc:Well-adapted_models}. In light of this it can be said that any manifold yields an example of an infinitesimal model of an affine space over $\mathbb{R}$.

  Although clones and their infinitesimal algebras are important for the discussion of properties of infinitesimal models of algebraic theories, in general, they are not convenient structures when one wants to work with infinitesimal models of a particular algebraic theory in practice. Indeed, having to work with the clone $\mathbb{T}$ is rather cumbersome. Fortunately, it is not necessary.
  
  Given a presentation $(\Sigma,E)$ of the algebraic theory $\mathbb{T}$ one can extend the signature and axioms of $\mathbb{T}$ to a new theory $\IIT$, such that $O_\mathbb{T}\iAlg\cong \IIT\Mod$; the latter denoting the category of $\IIT$-models and corresponding homomorphisms of $\Sigma$ structures. The new theory $\IIT$ is called the \emph{infinitesimalisation} of $\mathbb{T}$. It is obtained by 
  \begin{enumerate}[(1)]
    \item adding the theory of an i-structure,
    \item restricting the defining operations in $\Sigma$ to the i-structure,
    \item quantifying the defining equations over the i-structure,
    \item Adding a neighbourhood axiom for each defining operation $\Sigma$.
  \end{enumerate}
  More formally, and if one wishes to remain within the cartesian fragmet of first-order logic (cf. \cite[chap.~D1]{Johnstone:Elephant}), one needs to add relation symbols $A\la n\ra$, $n\in \IN$ to the signature, replace each function symbol with a functional relation (i.e. its graph) and then replace each $\Sigma$-term by a formula build from functional relations in conjunction with the respective i-structure before adding the i-structure and neighbourhood axioms. (For the technical details cf. \cite[def.~2.4.1]{Bar:thesis}).\footnote{The original infinitesimalisation construction given in \cite[def.~2.4.1]{Bar:thesis} does not remove the burden of proving the neighbourhood axiom for all derived operations over the signature $\Sigma$. For the simplified neighbourhood axioms mentioned in~(4) see \cite[def.~3, lem.~2 \& rem.~1]{Bar:2nd-order_i-groups}.} 

  Note that $\IIT$ is not algebraic anymore; it is a cartesian theory. Since the i-structure does not need to be defined equationally\footnote{The reason we do not make this restriction is due to examples like the formal manifolds in Synthetic Differential Geometry. Their i-structures are glued together from local models and are not globally defined by equations.}, $\IIT$ is also not essentially algebraic, in general.
 
  \begin{proposition}\label{prop:i-sation}
    Let $\mathbb{T}$ be an algebraic theory over a signature $\Sigma$ and $O_\mathbb{T}$ its clone as in proposition~\ref{prop:compare_syntactic_clone}. The categories $O_\mathbb{T}\iAlg$ and $\IIT\Mod$ are isomorphic as categories over $\Set$.
  \end{proposition}
  \begin{proof}
    It is not difficult to see that infinitesimal $O_\mathbb{T}$-algebras correspond to $\IIT$-models and vice versa. However, due to the formal logic involved the proof is rather lengthy and technical. The interested reader is referred to \cite[thm.~2.4.2]{Bar:thesis}, where this is proven for finite-limit categories as the base.
  \end{proof} 
  
  A pleasant consequence of the infinitesimalisation construction, respectively proposition~\ref{prop:i-sation}, is that the category $O\iAlg$ of infinitesimal $O$-algebras is locally finitely presentable:

  \begin{theorem}\label{thm:local_presentability}
    The category $O\iAlg$ is locally finitely presentable.  
  \end{theorem}
  \begin{proof}
    By proposition~\ref{prop:compare_syntactic_clone}, and since categories of infinitesimal algebras of isomorphic clones are isomorphic, we can assume w.l.o.g. that $O=O_\mathbb{T}$ for an algebraic theory $\mathbb{T}$.  Due to proposition~\ref{prop:i-sation} we have that $O_\mathbb{T}\iAlg\cong\IIT\Mod$. Since $\IIT$ is a cartesian theory and thus a (finitary) limit theory its category of models is locally finitely presentable \cite[thm.~5.9]{Adamek_Rosicky:LPAC}.  
  \end{proof}

  In particular, $O_\mathbb{T}\iAlg$ is complete, cocomplete, well-powered and well-copowered \cite[rem.~1.56]{Adamek_Rosicky:LPAC}. 

  \begin{theorem}\label{thm:lifting_of_lim_and_filtered_colim}\hspace{1ex}
    \begin{enumerate}[(i)]
      \item The forgetful functor $U:O\iAlg\to \Set$ lifts small limits uniquely; that is, for every small diagram $D:J\to O\iAlg$ and limiting cone $\lambda$ of $U\circ D$, there is a unique limiting cone $\mu$ of $D$ such that $U\mu=\lambda$. (See also \cite[def.~13.17]{Adamek_etal:ACC}.) 
      \item $U$ lifts filtered colimits uniquely.
    \end{enumerate}
  \end{theorem} 
  
  \begin{proof} 
    The infinitesimal $O$-algebras are constructed using products, monomorphisms and equalisers (for the defining equations) in $\Set$, which all commute with taking limits and filtered colimits, so $U$ lifts them. Since the functor $U$ is amnestic, i.e. an isomorphism $h$ is the identity morphism if $Uh$ is the identity map, all the limits and colimits $U$ lifts, it lifts uniquely. (See also \cite[prop.~13.21]{Adamek_etal:ACC}. Note that $U$ does neither reflect identities nor isomorphisms, in general, though.)

    More formally, we can use that $O\iAlg$ is equivalent to the category $Lex(C_{\IIT},\Set)$ of finite-limit preserving functors $C_{\IIT} \to \Set$ (for $\mathbb{T}=\mathbb{T}_O$), where $C_{\IIT}$ denotes the syntactic category of $\IIT$ \cite[thm.~D1.4.7]{Johnstone:Elephant}. The equivalence of categories commutes with the forgetful functors to $\Set$. We can thus consider $U$ to be the forgetful functor $\ev_A:Lex(C_{\IIT},\Set)\to\Set$, which is the evaluation at (the syntactical representation of) the (unique) sort $A$ in $C_{\IIT}$. Since limits and filtered colimits commute with finite limits, the limits and filtered colimits of finite-limit preserving functors are computed pointwise, which shows that $\ev_A$ lifts both.
  \end{proof}

  \begin{corollary}[Free i-$O$-algebras]
    $U$ has a left adjoint.
  \end{corollary}
  \begin{proof}
    $U$ preserves filtered colimits and small limits. By the adjoint functor theorem for locally presentable categories it has a left adjoint \cite[1.66]{Adamek_Rosicky:LPAC}. 
  \end{proof}

  As there are different i-structures on $A$, and therefore potentially more than one infinitesimal $O$-algebra structure, $U$ does not reflect limits, in general. For example, let $O$ be the clone of affine combinations over some field. Any affine space is an infinitesimal $O$-algebra for both the indiscrete and discrete i-structure, so $U$ does not reflect identities and thus does not reflect limits. This is a property infinitesimal models of algebraic theories share with topological spaces. In particular, as $U$ doesn't reflect isomorphisms, it fails to be monadic nor is it an essentially algebraic functor as defined in \cite[def.~23.1]{Adamek_etal:ACC}.
  
%%%%%%%%%%%%%%%%%%%%%%%%%%%%%%%%%%%%%%%%%%%%%%%%%%%%%%%%%%%%%%%%%%%%%%%%%%%%%%%%%%%%%%%%%%%%%%%%%%%%%

  \section{Gluing theorems for infinitesimal models}    

  We have seen that passing from the category of $\mathbb{T}$-models to infinitesimal models of an algebraic theory $\mathbb{T}$ the forgetful functor $U$ to the base category $\Set$ turns from being finitary monadic and algebraic \cite[def.~23.19]{Adamek_etal:ACC} to a functor that is neither monadic, nor essentially algebraic. ($U$ is not a topological functor \cite[def.~21.1]{Adamek_etal:ACC}, either.) 
  Although we loose the `algebraicity' over $\Set$ in the categorical sense, we retain many good categorical properties like local presentability and that limits and filtered colimits are computed from the underlying sets. In this section we study to which extent this also holds true for colimits.
  
  We begin by giving an explicit representation of the initial object in the category of infinitesimal models. 
  \begin{proposition}\label{prop:initial_object} Let $O$ be a clone.
    \begin{enumerate}[(i)]
      \item The set of constants $O(0)$ together with the maps $\ast_{(-)0}:O(n)\to O(0)$ is a (total) $O$-algebra.   
      \item The total $O$-algebra $O(0)$ is an initial object in $O\iAlg$.
    \end{enumerate}
  \end{proposition}
  \begin{proof}
    \begin{enumerate}[(i)]
      \item This is a consequence of the associativity (for $m=0$) and projection axioms of a clone $O$ given in definition~\ref{def:clone}.
      \item Let $A$ be an infinitesimal $O$-algebra. The neighbourhood axiom of definition~\ref{def:i-alg} shows (for $m=0$) that any $n$-tuple of constants in $O(0)$ lies in $A\la n\ra$; the map $\bullet_0:O(0)\to A\la 1\ra$ is thus an i-morphism. Since $A\la 0 \ra\cong 1$ the associativity for infinitesimal $O$-algebras shows $\bullet_0$ an i-$O$-homomorphism. The uniqueness of the i-$O$-homomorphism $(O(0),\ast_{(-)0}) \to (A,\bullet)$ is a consequence of $O$ satisfying the normalisation axiom. 
    \end{enumerate}  
  \end{proof}

  \begin{corollary}\label{cor:initial} The forgetful functor $U:O\iAlg\to \Set$ preserves the initial object if and only if the signature $\Sigma$ of $\mathbb{T}_O$ has no constants; i.e. $O(0)=\emptyset$. In that case $O\iAlg$ has a strict initial object and $U$ lifts and relfects it.
  \end{corollary}  

  Bearing in mind that arbitrary small colimits in a category can be constructed from the initial object and (small) wide pushouts the best result we could hope for is that $U$ lifts wide pushouts uniquely. It turns out that this is indeed the case provided the i-$O$-homomorphisms in the wide span \emph{reflect i-structure}.

  \begin{definition}
    Let $A$ and $B$ be i-structures. We say that an i-morphism $h:A\to B$ \textbf{reflects i-structure}, if it satisfies
    $$ \la h(x_1),\ldots,h(x_n)\ra \in B\la n\ra \implies\ \la x_1,\ldots,x_n\ra \in A\la n\ra$$ 
    for all $n\in\IN$. (The case $n=0$ is trivial.) An i-$O$-homomorphisms reflects i-structure, if the underlying i-morphism does.
  \end{definition}

  Any i-$O$-homomorphism from a total $O$-algebra is necessarily i-structure reflecting. However, not every i-$O$-homomorphism is i-structure reflecting. For example, the i-$O$-homomorphism from the discrete i-structure on an affine space $A$ to the indiscrete i-structure induced by the identity map we have encountered before does not reflect i-structure as long as $A$ is not the one-point set.  

  We are now ready to state and prove our first (and main) gluing theorem for infinitesimal models of an algebraic theory.

  \begin{theorem}\label{thm:lift_of_wide_pushouts}
    The forgetful functor $U:O\iAlg\to\Set$ lifts (small) wide pushouts of i-structure reflecting i-$O$-homomorphisms uniquely. Moreover, all the i-$O$-homomorphisms of the colimiting cocone reflect i-structure. 
  \end{theorem}
  \begin{proof}
    We will show that the assertion of the theorem holds true for (binary) pushouts first, and explain how the given construction generalises to the case of wide pushouts after. 
    
    \begin{asparaenum}[(1)]
      \item Consider the span of i-$O$-algebras, where $f$ and $g$ are i-structure reflecting i-$O$-homomorphisms.     
      $$
      \begin{tikzcd}
        & (A,\bullet^A)\ar[dl,"f"'] \ar[dr,"g"] & \\  
        (C,\bullet^C) & & (B,\bullet^B)
      \end{tikzcd}
      $$
      and a pushout of its $U$-image in $\Set$ 
      $$
      \begin{tikzcd}
        A\ar[d,"f"']\ar[r,"g"] & B\ar[d,"k"] \\
        C \ar[r,"h"] & Z
      \end{tikzcd}
      $$
      Recall that a pushout in $\Set$ can be constructed as the quotient of the coproduct $C\coprod B$ by an equivalence relation, which is generated by the relation $\set{(f(a),g(a))}{a\in A}$.
      $$
      \begin{tikzcd}
        A\ar[d,"f"']\ar[r,"g"] & B\ar[d,"i_B"] \\
        C \ar[r,"i_C"] & C\coprod_A B
      \end{tikzcd}
      $$ 
      The maps $i_C$ and $i_B$ are the coproduct inclusions composed with the quotient map. Two elements $i_C(x)$ and $i_B(y)$ are equal if and only if there are $a_1,\ldots,a_\ell$ in $A$ and a \emph{zig-zag} 
      $$
      \begin{tikzcd}[column sep = tiny]
        & a_1\ar[dl,"f"']\ar[dr,"g"]  & & a_2\ar[dl,"g"']\ar[dr,"f"] & & \ldots & & a_\ell\ar[dl,"f"']\ar[dr,"g"] &\\
      x=f(a_1) & & g(a_1) & & f(a_2) & \ldots & f(a_{\ell-1}) & & y=g(a_\ell) 
      \end{tikzcd}
      $$     
      and similarly for $i_X(x)$ and $i_X(x')$ as well as $i_Y(y)$ and $i_Y(y')$. (The zig-zags are just an explicit description of the transitive closure of (the reflexive symmetrisation of) the relation $\set{(f(a),g(a))}{a\in A}$.) Moreover, this characterisation by zig-zags is independent of the concrete choice of isomorphic representant of $Z$, $h$ and $k$.  
  
      \item To define an i-$O$-algebra structure on $Z$ we take the images of the i-$O$-algebra structures on $C$ and $B$ under $h$ and $k$, respectively. This makes $h$ and $k$ i-$O$-algebra homomorphisms that relfect i-structure by construction. 
      
      In particular, we define $Z\la n\ra$, $n\in\IN$ as the join of the images of $C\la n\ra$ and $B\la n\ra$ under $h^n$ and $k^n$, respectively. Since $h$ and $k$ are jointly epimorphic this defines an i-structure on $Z$. 
      
      For each $n\in\mathbb{N}$, $\sigma\in O(n)$ and $\la z_1,\ldots, z_n\ra\in Z\la n\ra$ we set $\sigma\bullet^Z_n (z_1,\ldots z_n)$ to be $h(\sigma\bullet^C_n(x_1,\ldots,x_n))$ or $k(\sigma\bullet^B_n(y_1,\ldots,y_n))$, for any $x_j$ or $y_j$ such that $z_j=h(x_j)$ or $z_j=k(y_j)$, respectively. We need to show that each $\bullet^Z_n$ yields a well-defined map $\bullet^Z_n:O(n)\times Z\la n\ra \to Z$. 
  
      Firstly, from the construction of the i-structure on $Z$ it follows that $\sigma\bullet^Z_n (z_1,\ldots z_n)$ is indeed defined for each $\la z_1,\ldots, z_n\ra\in Z\la n\ra$. In the case that $z_j=h(x_j)=k(y_j)$, there are zig-zags from $x_j$ to $y_j$ for each $1\leq j\leq n$. As zig-zags can be extended trivially, 
      % Always in steps of 2
      we can assume they are all of the same length $\ell$. From each zig-zag we take the first vertex $a^j_1\in A$. Due to $x_j=f(a^j_1)$, $\la x_1,\ldots,x_n \ra\in B\la n\ra$, and since $f$ reflects i-structure, $\la a^1_1,\ldots,a^n_1\ra \in A\la n\ra$. The morphisms $g$ and $f$ both preserve and reflect i-structure, so an easy induction over the length $\ell$ of the zig-zags shows $\la a^1_i,\ldots,a^n_i\ra \in A\la n\ra$ for each $1\leq i\leq \ell$. Since $f$ and $g$ are i-$O$-homomorphisms this yields a zig-zag    
      $$
      \begin{tikzcd}[column sep = 0.3em]
        & \sigma\bullet_n^A(\vec{a}_1) \ar[dl,"f"']\ar[dr,"g"]   & & \ldots & &\sigma\bullet^A_n(\vec{a}_{\ell})\ar[dl,"f"']\ar[dr,"g"] &\\
      \sigma\bullet^B_n (\vec{x}) & & \sigma\bullet^C_n g^n(\vec{a}_1)  &\ldots&  \sigma\bullet^B_n f^n(\vec{a}_{\ell}) & & \sigma\bullet^C_n (\vec{y})  
      \end{tikzcd}
      $$     
      and thus $\sigma\bullet^B_n (\vec{x}) = \sigma\bullet^C_n (\vec{y})$. The other cases, where $z_j=h(x_j)=h(x'_j)$, or $k(y_j)=k(y'_j)$ for all $j$ can be treated the same way. 

      It is now straight-forward to verify that $(Z,\bullet^Z)$ is an i-$O$-algebra. Indeed, any $\la z_1,\ldots,z_n\ra$ can be represented as an $h^n$-image of some $\la x_1,\ldots,x_n\ra \in B\la n\ra$ or $k^n$-image of some $\la y_1,\ldots,y_n\ra \in C\la n\ra$. Since all the axioms are equations holding true in $(B,\bullet_B)$ and $(C,\bullet^C)$, they also satisfied by $(Z,\bullet^Z)$. Moreover, $h$ and $k$ are i-structure reflecting i-$O$-homomorphisms by construction.
      
      In the same vein it follows that $(Z,\bullet^Z)$ together with $h$ and $k$ is a pushout in $O\iAlg$: For an i-$O$-algebra $(W,\bullet^W)$ and i-$O$-homomorphisms $r$ and $s$ the unique map $t$ in the commutative diagram
      $$
      \begin{tikzcd}
        A\ar[d,"f"']\ar[r,"g"] & B\ar[d,"k"]\ar[ddr,"s",bend left] &\\
        C\ar[drr,"r"',bend right] \ar[r,"h"] & Z\ar[dr,dashed,"t"] & \\
        & & W
      \end{tikzcd}
      $$
      lifts to an i-$O$-homomorphism $t:(Z,\bullet^Z)\to (W,\bullet^W)$ by virtue of the construction of the i-$O$-algebra structure on $Z$.

      \item The general case of (small) wide pushouts follows from generalising the construction in (2) to an arbitrary (small) set of i-$O$-homomorphisms $I$ with common domain $(A,\bullet^A)$. Firstly, if we denote the codomain of each map $f\in I$ by $A_f$, then the wide pushout can be constructed as the quotient of the coproduct $\coprod_{f\in I} A_f$ by the equivalence relation generated by the relation $\set{(f(a),g(a))}{a\in A,\ f,g\in I}$. As in the binary case, the transitive closure of this relation can be represented by zig-zags for all pairs $(f,g)\in I^2$, and this representation yields a characterisation that is independent of the chosen construction of the wide pushout.  
      
      The i-$O$-algebra structure on the wide pushout $Z$ with maps $i_f:A_f\to Z$ is constructed as in the binary case: we take the images of the i-$O$-algebra structures under all $f\in I$. By applying the argument given in (2) to all pairs $(f,g)\in I^2$, $f\neq g$, we see this to be well-defined. As in (2), the $i_f$ are i-structure reflecting i-$O$-homomorphisms and satisfy the universal property of a wide pushout in $O\iAlg$ by construction. The lift of the wide pushout is necessarily unique due to $U$ being amnestic (see the proof of theorem~\ref{thm:lifting_of_lim_and_filtered_colim}).
    \end{asparaenum} 
  \end{proof}

  As the initial object $O(0)$ is a total $O$-algebra the unique i-$O$-homomorphisms reflect i-structure. Theorem~\ref{thm:lift_of_wide_pushouts} thus provides us with an explicit construction of coproducts of i-$O$-algebras as (wide) pushouts with vertex in $O(0)$. 
  
  By considering the coslice category $O(0)/\Set$ of maps $f$ with $\dom f=O(0)$ and morphisms $f\to g$ maps $h:X\to Y$ with $hf=g$
  as the base category we can phrase this result as follows:
  
  \begin{corollary}\label{cor:mapping_coproducts}
    The forgetful functor $U:O\iAlg\to O(0)/\Set$ lifts (small) coproducts uniquely. In particular, if $O(0)=\emptyset$, then coproducts in $O\iAlg$ are disjoint unions.
  \end{corollary}
  
  Recall that a coequalizer $q$ of $f,g:A\rightrightarrows B$ can be written as a pushout
  $$
  \begin{tikzcd}
    A\coprod B\ar[d,"{(f,1_B)}"'] \ar[r,"{(g,1_B)}"] & B\ar[d,"q"] \\
    B\ar[r,"q"] & Z
  \end{tikzcd}
  $$
  However, the folding maps $(f,1_B)$ and $(g,1_B)$ are not necessarily i-structure reflecting, even if $f$ and $g$ are. To guarantee this property the images of $f$ and $g$ have to be \emph{infinitesimally closed} subsets; i.e. if $\langle f(x), y\rangle\in B\langle 2\rangle$ then $y$ has to lie in the image of $f$ (and the same for $g$). Indeed, $(f,1_B)$ is i-structure reflecting if and only if $f$ is i-structure reflecting and has an infinitesimally closed image in $B$. 
  Another consequence of theorem~\ref{thm:lift_of_wide_pushouts} is therefore:

  \begin{corollary}\label{cor:lifting_coequalizers}
    The forgetful functor $U:O\iAlg\to \Set$ lifts the coequalizer of a (small) family of parallel i-structure reflecting i-$O$-homomorphisms uniquely, if each of the maps has an infinitesimally closed image. In this case the coequalizer reflects i-structure. 
  \end{corollary}
  
  Since (small) colimits can be constructed systematically from (small) coproducts and a coequalizer (cf. the dual of \cite[thm.~V.2.1]{MacLane:CWM}), one might hope to be able conclude that $U$ lifts (small) colimits of i-structure reflecting i-$O$-homomorphisms with infinitesimally closed images from theorem~\ref{thm:lift_of_wide_pushouts}. However, it turns out that the representation of the colimit given by the dual of \cite[thm.~V.2.1]{MacLane:CWM} is unsuitable, since the maps between the coproducts will not reflect i-structure. The problem lies with the folding maps from a coproduct used in this construction. Consider, for example the codiagonal folding map $\delta=(1_A,1_A)$ for an i-$O$-algebra $A$
  $$
  \begin{tikzcd}
    A \arrow[d, tail] \arrow[rrd, "1_A"]   &  &   \\
    A\coprod A \arrow[rr, "\delta" near start]             &  & A \\
    A \arrow[u, tail] \arrow[rru, "1_A"'] &  &  
  \end{tikzcd}
  $$
  It is i-structure reflecting if and only if $A$ is a subalgebra of $O(0)$. 
  We shall provide an example of a pair of two i-structure reflecting i-$O$-homomorphisms for which $U$ does not lift the coequalizer later, so corollary~\ref{cor:lifting_coequalizers} is already the best possible general result. 
  
  In our second gluing theorem we show that coequalizers in $O\iAlg$ can be constructed as in $O\Alg$, in general: namely, as quotients of congruences\footnote{Recall that a congruence is an equivalence relation compatible with the structure. In the case of $O\iAlg$ a congruence on $(A,\bullet^A)$ amounts to two i-$O$-homomorphisms $(p_1,p_2): (R,\bullet^R)\rightrightarrows (A,\bullet^A)$ that are jointly monomorphic, such that the induced $UR\monto UA \times UA$ is an equivalence relation in $\Set$.}. 
  
  \begin{theorem}\label{thm:lifts_of_quotients}
    Let $
    \begin{tikzcd}
      (R,\bullet^R)\ar[r,shift left,"p_1"]\ar[r,shift right,"p_2"'] & (A,\bullet^A)
    \end{tikzcd}
    $
    be a congruence in $O\iAlg$. $U$ lifts the quotient of $R$ uniquely if $p_1$ and $p_2$ jointly reflect the i-structure; that is for all $n\in\mathbb{N}$ and $y_1,\ldots,y_n\in R$
    $$ 
      \la p_1(y_1),\ldots,p_1(y_n)\ra \in A\la n \ra \wedge \la p_2(y_1),\ldots,p_2(y_n)\ra \in A\la n \ra\implies \la y_1,\ldots,y_n\ra\in R\la n\ra
    $$
    In this case the lift of the quotient map reflects i-structure.
  \end{theorem} 
  \begin{proof}
    $R$ is also an equivalence relation in $\Set$. Let $q$ be its quotient. 
    $$
    \begin{tikzcd}
       R\ar[r,shift left,"p_1"]\ar[r,shift right,"p_2"'] & A \ar[r,"q"] & Z
    \end{tikzcd}
    $$
    Firstly, we show that $p_1$ and $p_2$ jointly reflecting i-structure is sufficient for the existence of a $U$-lift of the coequalizer, which is then necessarily unique, for $U$ is amnestic. As in the proof of theorem~\ref{thm:lift_of_wide_pushouts} the i-$O$-structure on $Z$ is constructed by taking the respective images of the i-$O$-algebra structure under $q$. For $n\in \mathbb{N}$ the i-structure $Z\la n\ra$ is thus defined as $q^n(A\la n\ra)$ and the $O$-action $\bullet^Z$ is defined as 
    $$
      \sigma\bullet^Z(z_1,\ldots, z_n) = q(\sigma\bullet^A(x_1,\ldots,x_n))
    $$
    for any $\sigma\in O(n)$ and $\la x_1,\ldots,x_n\ra\in A\la n\ra$ such that $q(x_j)=z_j$. 
            
    Since $A\la -\ra$ is an i-structure and $q$ an epimorphism, it is clear that $Z\la -\ra$ defines an i-structure for which $q$ is an i-structure reflecting i-morphism. Moreover, $\sigma\bullet^Z_n(z_1,\ldots, z_n)$ is defined for every $\la z_1,\ldots,z_n\ra\in Z\la n\ra$. We need to show that $\sigma\bullet_n^Z(z_1,\ldots,z_n)$ is well-defined. 

    Let $x'_j\in A$ be such that $q(x'_j)=q(x_j)$ and $\la x'_1,\ldots,x'_n \ra \in A\la n\ra$. Since $\Set$ is an effective regular category\footnote{Recall that a category is called regular, if it has finite limits and regular epimorphisms (i.e. epimorphisms that are coequalizers) are stable under pullback. In an effective regular category every congruence is a kernel pair.}, 
    $R$ is the kernel pair of $q$, and $R(x_j,x_j')$ holds; i.e. $x_j$ and $x_j'$ are $R$-equivalent. There are thus (uniquely determined) $y_j\in R$ such that $p_1(y_j)=x_j$ and $p_2(y_j)=x_j'$. Moreover, since both $\la x_1,\ldots,x_n\ra$ and $\la x_1',\ldots,x_n'\ra$ lie in $A\la n\ra$ and $p_1$ and $p_2$ jointly reflect i-structure, we have $\la y_1,\ldots,y_n\ra\in R\la n\ra$. As $p_1$ and $p_2$ are i-$O$-homomorphisms we find
    \begin{align*}
      \sigma\bullet_n^A(x_1,\ldots,x_n) &= p_1(\sigma\bullet_n^R(y_1,\ldots,y_n)) \\
      \sigma\bullet_n^A(x'_1,\ldots,x'_n) &= p_2(\sigma\bullet_n^R(y_1,\ldots,y_n))
    \end{align*}
    and thus $q(\sigma\bullet_n^A(x_1,\ldots,x_n))=q(\sigma\bullet_n^A(x_1',\ldots,x_n'))$, for $q$ is the coequalizer of $p_1$ and $p_2$.
     
    As in the proof of theorem~\ref{thm:lift_of_wide_pushouts} the axioms of an i-$O$-algebra for $(Z,\bullet^Z)$ follow easily from the definition of the i-$O$-algebra structure on $q$ and the fact that they hold for $(A,\bullet^A)$. Clearly, $q$ lifts to an i-$O$-homomorphism $(A,\bullet^A)\to (Z,\bullet^Z)$ and the universal property of $q$ as the coequalizer of $p_1$ and $p_2$ in $O\iAlg$ follows easily as well.        
  \end{proof}
  
  Note that the converse of theorem~\ref{thm:lifts_of_quotients} does not hold, and requires the additional assumption that the equivalence relation $(p_1,p_2):(R,\bullet^R)\rightrightarrows (A,\bullet^A)$ is a kernel pair in $O\iAlg$. Indeed, let $A$ be an affine space with more than one point considered as a total i-$O$-algebra $(A,\bullet^A)$ for the clone $O$ of affine combinations. Consider the discrete i-structure on the set $A\times A$ making it into an i-$O$-algebra $(A\times A,\delta)$. The pair of projections $(\pr_1,\pr_2):(A\times A, \delta)\rightrightarrows (A,\bullet^A)$ is a congruence, for which $p_1$ and $p_2$ do not jointly reflect i-structure, but its quotient $UA\to 1$ in $\Set$ has a unique lift\footnote{Note that the "only if" part stated in theorem~2.6.20 in \cite{Bar:thesis} is incorrect for exactly this reason. Also, the statement that $U$ reflects regular epis in \cite[cor.~2.6.21]{Bar:thesis} is false; the counter example being provided in remark~2.6.22(a) ibid.}.

  \begin{corollary}\label{cor:characterising_kernels}
    A congruence 
    $
    \begin{tikzcd}
      (R,\bullet^R)\ar[r,shift left,"p_1"]\ar[r,shift right,"p_2"'] & (A,\bullet^A)
    \end{tikzcd}
    $ in $O\iAlg$ is a kernel pair if and only if $p_1$ and $p_2$ jointly reflect i-structure.
  \end{corollary}
  \begin{proof}
    \begin{asparaenum}[(1)]
      \item Let $(p_1,p_2): (R,\bullet^R)\rightrightarrows (A,\bullet^A)$ be a kernel pair. Since $U$ preserves but also uniquely lifts kernel pairs, the i-structure on $R$ is the restriction of the product i-structure on $A\times A$ along the monic $(p_1,p_2):R \monto A\times A$; but this shows $p_1$ and $p_2$ to jointly reflect i-structure.
      
      A more formal argument can be given as follows: The assertion that $p_1$ and $p_2$ jointly reflect i-structure is equivalent to saying that the commutative squares
      $$
        \begin{tikzcd}
          R\la n\ra \ar[d,rightarrowtail]\ar[r,rightarrowtail,"{(p_1^n,p_2^n)}"] & A\la n\ra\times A\la n\ra \ar[d,rightarrowtail] \\
          R^n \ar[r,rightarrowtail,"{(p_1^n,p_2^n)}"] & A^n\times A^n
        \end{tikzcd}		
      $$
      are pullbacks (in $\Set$) for $n\in\IN$. This follows from a diagram chase 
      $$
      \begin{tikzcd}
        R\la n\ra\ar[d,rightarrowtail]\ar[r,shift left,"{p^n_1}"]\ar[r,shift right,"{p^n_2}"'] & A\la n\ra\ar[d,rightarrowtail] \ar[r, twoheadrightarrow,"q^n"] & Z\la n\ra \ar[d,rightarrowtail]\\
        R^n\ar[r,shift left,"{p^n_1}"]\ar[r,shift right,"{p^n_2}"'] & A^n \ar[r, twoheadrightarrow,"q^n"] & Z^n
      \end{tikzcd}
      $$ 
      where $q$ is the coequalizer of $(p_1,p_2)$ in $\Set$, using the fact that the top and bottom diagrams are kernel pairs in $\Set$. (The i-structure on $Z$ is the $q$-image of the i-structure on $A$. It can be defined irrespective of whether $Z$ can be made into an i-$O$-algebra.)

      \item Conversely, from the proof of the second gluing theorem~\ref{thm:lifts_of_quotients} we can see that $(R,\bullet^R)$ is the kernel pair of the lift of the quotient $q:UA\to Z$. Indeed, since $U$ lifts the kernel pair $(Up_1, Up_2): UR\rightrightarrows UA$ of $Uq:UA\to UZ$ uniquely (we shall denote this lift by $(R',\bullet')$), there is a unique i-$O$-homomorphism $h: (R,\bullet^R)\to (R',\bullet')$ satisfying $p_j\circ h = p_j$ and $Uh$ being the identity map on $R=R\la 1\ra =R'\la 1\ra$. However, since $p_1$ and $p_2$ jointly reflect i-structure (for both $R$ and $R'$), they have to be equal; namely the restriction of the product i-structure on $A\times A$ along the monic $(p_1,p_2): R \monto A \times A$ as argued in (1). This shows that $h$ is the identity morphism in $O\iAlg$. 
    \end{asparaenum} 
  \end{proof}

  Since $U$ preserves the equivalence relation $(R,\bullet^R)\rightrightarrows (A,\bullet^A)$ and every equivalence relation in $\Set$ is a kernel pair, part~(2) of the proof of the preceding corollary shows that the only difference between $(R,\bullet^R)$ and the kernel pair $(R',\bullet')$ is the i-structure (since $Uh$ has to be the identity map); it is only if the i-structure is `too small' relative to the codomain that prevents an equivalence relation in $O\iAlg$ to be a kernel pair. 
  
  \begin{corollary}\label{cor:OiAlg_regular}
    $O\iAlg$ is a regular category (but not effective regular) and $U$ is a regular functor, i.e. preserves finite limits and regular epimorphisms.
  \end{corollary}
  \begin{proof}
    \begin{asparaenum}[(1)]
      \item Recall that a regular epimorphism $e$ is a coequalizer of its kernel pair. By the preceding corollary~\ref{cor:characterising_kernels} each kernel pair jointly reflects i-structure. By theorem~\ref{thm:lifts_of_quotients}, $U$ lifts its $\Set$-quotient uniquely. The lifted quotient is thus isomorphic to the coequalizer $e$ in $O\iAlg$ and $Ue$ is thus a coequalizer in $\Set$.
      
      \item We have seen in the preceding discussion that $O\iAlg$ has equivalence relations that are not kernel pairs, so it cannot be an effective regular category. To show it a regular category, since $O\iAlg$ is complete and cocomplete, it is sufficient to show that regular epimorphisms are stable under pullbacks. 
    
      Now consider a pullback $f^*e$ of a regular epimorphism $e$ along $f$ in $O\iAlg$: 
      
      \begin{equation}
          \begin{tikzcd}
            C \ar[d,rightarrow,"f^*e"']\ar[r,rightarrow,"e^*f"] & A \arrow[d,twoheadrightarrow,"e"] \\
            B \ar[r,rightarrow,"f"] & Z
            \label{diag:pb_of_reg_epi}
          \end{tikzcd}		
      \end{equation}  
      $U$ preserves pullbacks and regular epis, so the $U$-image of the diagram is a pullback diagram of the regular epi $Ue$ in $\Set$. Since $\Set$ is a regular category $U(f^*e)$ is a regular epi in $\Set$ and the coequalizer of its kernel pair. By theorems~\ref{thm:lifting_of_lim_and_filtered_colim} and \ref{thm:lifts_of_quotients}, $U$ lifts kernel pairs and thus their quotients uniquely. Applying the argument given in part~(2) of the proof of the preceding corollary~\ref{cor:characterising_kernels} to the lifted coequalizer here, we get that the lifted coequalizer agrees with $f^*e$ if and only if the i-structure on $B$ is the image of the i-structure on $C$ under $f^*e$.
      
      This follows from the following diagram chase utilising the construction of the i-structures on $Z$ and $C$: Let $n>1$ and $\la b_1,\ldots,b_n\ra\in B\la n\ra$. Since $e$ reflects i-structure, there is $\la a_1,\ldots,a_n\ra\in A\la n\ra$ with $e(a_j)=f(b_j)$ for all $1\leq j\leq n$. As the above diagram~(\ref{diag:pb_of_reg_epi}) is a pullback there exist $c_j\in C$ such that $f^*e(c_j)=b_j$ and $e^*f(c_j)=a_j$. Moreover, since $f^*e$ and $e^*f$ jointly reflect i-structure, we have $\la c_1,\ldots,c_n\ra\in C\la n\ra$. This shows $f^*e$ a regular epi.
    \end{asparaenum}      
  \end{proof}

  \begin{corollary}\label{cor:inclusion_of_O-Alg}
    The (full) inclusion $O\Alg\into O\iAlg$ mapping $O$-algebras to total i-$O$-algebras preserves (small) limits and coequalizers, but does not preserve coproducts, in general.
  \end{corollary}
  \begin{proof}
    \begin{asparaenum}[(1)]
      \item (Small) products of total i-$O$-algebras are total and the same holds true for equalizers. By \cite[thm.~V.2.1]{MacLane:CWM} (small) limits of diagrams of total i-$O$-algebras are total. As limits of both $O$-algebras and i-$O$-algebras are constructed from the underlying limits in $\Set$ (each of the forgetful functors lifts limits uniquely), this is equivalent to saying that the inclusion functor preserves (small) limits.

      \item A similar argument shows that coequalizers are preserved. Indeed, as seen in the proof of the preceding corollary~\ref{cor:OiAlg_regular} the coequalizer of a parallel pair of i-$O$-homomorphisms $f,g:A\rightrightarrows B$ is the $U$-lift of the quotient $q$ of its kernel pair. Essentially by the same argument this is also true for coequalizers of $O$-algebras, as the respective forgetful functor lifts $\Set$-quotients of congruences uniquely. From the construction of the i-structure in theorem~\ref{thm:lifts_of_quotients} we see that the quotient i-$O$-algebra will be total, if $B$ is total. In this case the kernel pair of $q$ is also total. This shows that the inclusion preserves coqualizers.  

      \item Finally, we note that the underlying set of a coproduct of two nontrivial abelian groups is a binary product of sets, which is different from the wedge sum of the underlying pointed sets, i.e. the underlying set of their coproduct as infinitesimal abelian groups. In the case of $O(0)=\emptyset$, we note that the coproduct of two affine lines is a three dimensional affine spaces and thus different from the disjoint union of the underlying sets.
    \end{asparaenum}
  \end{proof}

  With corollary~\ref{cor:inclusion_of_O-Alg} we see that the forgetful functor $U:O\iAlg\to\Set$ cannot lift coequalizers of i-structure reflecting i-$O$-homomorphisms, in general, since the forgetful functor $O\Alg\to\Set$ does not lift coequalizers, in general. Indeed, every i-$O$-homomorphism between total i-$O$-algebras is i-structure reflecting. However, the coequalizer of the the constant-zero map and $z\mapsto 2z$ of the integers $\IZ$ is the canonical projection onto $\IZ/2\IZ$ in the category of abelian groups; but in $\Set$ it is the map to the set $2\IZ+1\cup\{0\}$ that collapses even integers to $0$. 
  
  A similar counterexample can be given for affine spaces over $\IR$.\footnote{This shows the statement of \cite[thm.~2.6.19]{Bar:thesis} to be false.}  
  Consider $\IR^2$ as an affine space over itself. The constant-zero map and $x\mapsto (x,0)$ are affine. In the category of affine spaces the coequalizer is the projection to $\IR\times\{0\}$, but in the category $\Set$ it is the map to the wedge sum of the two pointed half planes $(0,\infty)\times\IR\cup\{(0,0)\}$, $(-\infty,0)\times\IR\cup\{(0,0)\}$ joined at $(0,0)$ that collapses the line $\IR\times\{0\}$ to $(0,0)$. The unique map from this wedge sum to $\IR\times\{0\}$ compatible with the quotient maps is surjective but not injective.

  Our third and last gluing theorem concerns the gluing of local models, which is of importance in the application of infinitesimal algebra in Synthetic Differential Geometry\footnote{In line with the rest of the paper we give the definition for $\Set$ as a base category here. However, $\Set$ can be replaced with a Grothendieck topos $\calS$ and the proof of the gluing theorem directly transfers to $\calS$ as well.}.

  Let $M$ be a set. A poset $\calA\subset \Sub(M)$ together with a lift of $\calA$ along the forgetful functor $U:O\iAlg\to \Set$ is an \emph{infinitesimal $O$-atlas for $M$}, if 
  \begin{enumerate}[(1)]
    \item $\calA$ is stable under finite meets
    \item Each i-$O$-monomorphism in the $U$-lift is i-structure reflecting
    \item The union of $\calA$ is $M$.
  \end{enumerate}

  \begin{theorem}\label{thm:gluing_local_models}
    The forgetful functor $U:O\iAlg\to\Set$ lifts the union of any infinitesimal $O$-atlas $\calA$ uniquely.  
  \end{theorem}
  \begin{proof}
    Any union is the filtered colimit of finite unions of subobjects in $\calA$. As $U$ lifts filtered colimits uniquely by theorem~\ref{thm:lifting_of_lim_and_filtered_colim}, we only need to consider the case of finite unions.

    The finite join $\bigcup_{i=1}^n UV_i$ of the subobjects $UV_i\in \calA$ is the wide pushout
    $$
    \begin{tikzcd}
                     & \bigcap_{i=1}^n UV_i \arrow[ld, tail] \arrow[rd, tail] \arrow[d, tail] &                      \\
    UV_1 \arrow[rd, tail] & \cdots \arrow[d, tail]                                                & UV_n \arrow[ld, tail] \\
                     & \bigcup_{i=1}^n UV_i                                                   &                     
    \end{tikzcd}
    $$
    Since $\calA$ is closed under finite meets there is an i-$O$-algebra $W$ with $UW = \bigcap_{i=1}^n UV_i$ and all the monomorphisms are $U$-images of i-structure reflecting i-$O$-monomorphisms $W\monto V_i$. By our first gluing theorem~\ref{thm:lift_of_wide_pushouts} the functor $U$ lifts the union $\bigcup_{i=1}^n UV_i$ uniquely and all the i-$O$ monomorphisms reflect i-structure. Iterating this argument we can extend $\calA$ to a filtered poset by adding all the finite unions.     
  \end{proof}
  Although we formulated the gluing theorem for general algebraic theories, it seems most relevant if the theory has no constants, i.e. $O(0)=\emptyset$. Indeed, in the presence of constants each point of an i-$O$-algebra is infinitesimally close to every constant, so speaking of local models seems strange. 
  In the case of a theory without constants the coproducts in $O\iAlg$ are the disjoint unions in $\Set$. This is the case for the theory of affine spaces, for example, and makes gluing of infinitesimal models of affine spaces convenient.

%%%%%%%%%%%%%%%%%%%%%%%%%%%%%%%%%%%%%%%%%%%%%%%%%%%%%%%%%%%%%%%%%%%%%%%%%%%%%%%%%%%%%%%%%%%%%%%%%%%%%
  \section{Conclusion}

  Extending the category of models of an algebraic theory $\IT$ to the category of infinitesimal models helps us retain many good categorical properties like being locally finitely presentable as well as computing limits and filtered colimits from the underlying sets. Coproducts reduce to wedge sums over the set of constants, respectively to disjoint unions, if the theory does not require constants. Which further colimits can be computed from the underlying sets is determined by the infinitesimal structure, respectively the morphisms. There are mild conditions for lifting pushouts and the gluing of local models, but coequalizers are computed from quotients of congruences like in the category of models of $\IT$, in general; it is the infinitesimal structure that determines whether the infinitesimal models are closer to geometry, or to algebra. This shows infinitesimal models of algebraic theories a convenient interpolation between these two (meta-)concepts. 
  
  So far, infinitesimal algebra has been solely applied in Synthetic Differential Geometry revealing interesting relationships between the algebra and well-known differential geometric concepts and constructions at the infinitesimal level \cite{Kock:Synthetic_Geometry_Manifolds, Kock:Levi-Civita}, \cite[chap.~3]{Bar:thesis}, \cite{Kock:flatness, Bar:second_order_affine_structures}; this is subject to ongoing research. However, We hope that presenting it as a construct independent of its geometric roots, it can find use cases and applications in other fields, not necessarily related to geometry.
  
  The construction of infinitesimalisation is not restricted to algebraic theories and can be extended to any first-order theory. This leads to categories of infinitesimal models for each fragment of first-order logic. Whether and how this changes the properties of the categories of models, as well as applications of these constructions remain open questions subject to future research.  

%%%%%%%%%%%%%%%%%%%%%%%%%%%%%%%%%%%%%%%%%%%%%%%%%%%%%%%%%%%%%%%%%%%%%%%%%%%%%%%%%%%%%%%%%%%%%%%%%%%%%%

  \bibliography{references_with_links}
\end{document}